\documentclass[a4paper,11pt]{amsart}

\usepackage{amsmath,amssymb,amsfonts,amsthm,mathrsfs}
\numberwithin{equation}{section}
\usepackage[pdfusetitle]{hyperref}
\usepackage{bookmark}
\usepackage{verbatim}
\usepackage{tikz}

\hypersetup{
  colorlinks=true,
  linkcolor=blue,
  citecolor=blue,
  urlcolor=blue
}

\newtheorem{theorem}{Theorem}[section]
\newtheorem{lemma}[theorem]{Lemma}

\theoremstyle{definition}
\newtheorem{definition}[theorem]{Definition}

\newtheorem{remark}[theorem]{Remark}

\newcommand{\R}{\mathbb{R}}

\title[Hausdorff Measure Estimates for Free Boundaries]{$(n-1)$-Dimensional Hausdorff Measure Estimates for Free Boundaries of Infinity Laplace Equations}

\author{Yantian Chen}
\address[Yantian Chen]{School of Mathematics, Southeast University,
Nanjing 211189, P.R. China}
\email{213230018@seu.edu.cn}

\author{Feida Jiang$^*$}
\address[Feida Jiang]{School of Mathematics and Shing-Tung Yau Center of Southeast University,
Southeast University, Nanjing 211189, P.R. China;
Shanghai Institute for Mathematics and Interdisciplinary Sciences,
Shanghai 200433, P.R. China}
\email{jiangfeida@seu.edu.cn}

\thanks{$^*$Corresponding author.}

\hypersetup{pdfauthor={Yantian Chen and Feida Jiang}}

\subjclass[2020]{Primary 35J94; Secondary 35R35, 35B65, 28A78}

\keywords{Infinity Laplacian, sharp growth estimates, free boundary, Hausdorff measure estimates}

\date{}

\newif\ifreferencesused
\AddToHook{cmd/cite/before}{\global\referencesusedtrue}
\AddToHook{cmd/nocite/before}{\global\referencesusedtrue}

\begin{document}

\begin{abstract}
We investigate the geometric properties of free boundaries associated with nonnegative viscosity solutions to the infinity Laplace equation with strong absorption
\[
\Delta_\infty u=\lambda u_+^m,
\qquad \lambda>0,\quad 0<m<3.
\]
We first establish a sharp distance-to-free-boundary regularity estimate with an explicit optimal coefficient, which coincides with the constant in the non-degeneracy estimate. This exact matching yields a uniform interior tangent ball property at every sufficiently small scale for every free boundary point. As consequences, we obtain a sharp lower density bound of \(1/2\) for the positive phase, a locally finite covering of the free boundary by Lipschitz graphs, and quantitative local estimates for its \((n-1)\)-dimensional Hausdorff measure, thereby giving an affirmative answer to the question raised by Ara\'ujo, Leit\~ao, and Teixeira (J. Funct. Anal. {\textbf 270} (2016), no. 6, 2249–2267). The underlying barrier mechanism also suggests possible extensions to more general equations.
\end{abstract}

\maketitle

\section{Introduction and main results}
In this paper, we consider the geometric properties of free boundaries of the equation involving the infinity Laplacian operator, 
\begin{equation}
\Delta_\infty u
:= \sum_{i,j} \partial_i u\,\partial_{ij}u\,\partial_j u
= (Du)^T D^2u\,Du.
\end{equation}

More precisely, we consider nonnegative viscosity solutions to
\begin{equation}
\label{eq:dirichlet}
\begin{cases}
\Delta_\infty u=\lambda u_+^m & \text{in }\Omega,\\
u=\phi & \text{on }\partial\Omega,
\end{cases}
\end{equation}
where $\lambda>0$, $0<m<3$, and $\phi\in C(\partial\Omega)$ is nonnegative.
The case $m=0$ is related to the infinity-obstacle problem. The constant $\lambda>0$ here is called the Thiele modulus, which adjusts the ratio of reaction rate to diffusion-convection rate. 

The infinity Laplacian originates from Aronsson's study of absolutely minimizing Lipschitz extensions. Jensen then established the viscosity framework and uniqueness for the corresponding homogeneous Dirichlet problem. The regularity theory for infinity harmonic functions has been developed, featuring the work of Evans, Savin and Smart.

In the more subtle inhomogeneous case, Lu and Wang studied the existence and uniqueness of viscosity solutions provided the inhomogeneous term $f$ has a strict fixed sign, i.e. either $\inf f > 0$ or else $\sup f < 0$ in \cite{LuWang2008}.

Over the past decades, reaction-diffusion models governed by the infinity Laplacian have been extensively studied. In \cite{AraujoLeitaoTeixeira2016}, the equation featuring m-strong absorption. A typical feature of this kind of equations is a phenomenon of plateaus, or dead cores, where a nonnegative solution vanishes identically. Much attention has been paid to the properties of the solution near free boundary since it presents an improved regularity and non-degeneracy property. 

Also, the geometric properties of free boundary itself, defined by
\[
\Gamma(u):=\partial\{u>0\}\cap\Omega.
\]
is also of considerable interest.

In this direction, Ara\'ujo, Leit\~ao and Teixeira~\cite{AraujoLeitaoTeixeira2016} completed a systematic study of equation
\begin{equation}
  \Delta_\infty u=\lambda u_+^m.
\end{equation}

They established optimal growth estimates and non-degeneracy property near the free boundary, with the intrinsic growth exponent ${4}/(3-m)$.

More precisely, the results provided upper growth estimates near critical points and lower growth estimates at the closure of positive phase, both of $r^{\frac{4}{3-m}}$. Based on the estimates of the solution, they claimed that the free boundary is porous and for some universal $0<\sigma_0\le1$, 
\[
\mathcal H^{n-\sigma_0}_{\mathrm{loc}}\bigl(\Gamma(u)\cap K\bigr)<\infty
\]
See \cite[Corollary~6.2]{AraujoLeitaoTeixeira2016}.
In \cite[Remark~6.3]{AraujoLeitaoTeixeira2016}, they asked whether $\sigma_0$ can reach $1$ and leave it as an open problem.

Subsequent work extended the theory in several directions.
Diehl and Teymurazyan~\cite{DiehlTeymurazyan2020} obtained sharp regularity estimates near the free boundary for more general reaction terms and proved porosity provided suitable conditions of the right hand side term.
Lin and Liu~\cite{LinLiu2022} studied the Dirichlet problem with a
spatially varying absorption coefficient $\lambda(x)$ satisfying
$\inf_{\Omega}\lambda(x)>0$ and a nonpositive source term.

For H\'enon-type equations, the absorption coefficient can vanish. The model weight $\lambda(x)=|x|^\alpha$, with $\alpha>0$, vanishes at the origin. Related weighted equations were considered by Bezerra J\'unior, da Silva, Nascimento, and S\'a~\cite{BezerraDaSilvaNascimentoSa2024}.

In our earlier work~\cite{ChenJiang2026}, we distinguished free boundary points $x_0$ with $\lambda(x_0)>0$ from those with $\lambda(x_0)=0$. We established a degeneracy-set comparison principle and obtained sharp growth estimates at degenerate free boundary points with suitable inhomogeneous term.

In the present paper, we return to the constant-coefficient problem \eqref{eq:dirichlet} and address the question raised in \cite[Remark~6.3]{AraujoLeitaoTeixeira2016} concerning the local finiteness of the $(n-1)$-dimensional Hausdorff measure of the free boundary. The key observation of this work is that the upper and lower growth estimates can be matched with the same sharp coefficient
\[
\tau:=\tau(\lambda,m)
=
\left(
\frac{\lambda(3-m)^4}{64(1+m)}
\right)^{\frac{1}{3-m}}.
\]

We establish the sharp upper growth estimate
by a direct comparison argument with radial
barriers. Combined with the optimal-constant
non-degeneracy estimate, this yields a uniform
interior tangent ball property at every free
boundary point.

This geometric property allows us to cover
the free boundary locally by finitely many
Lipschitz graphs, leading to quantitative
Hausdorff measure estimates.

Our main result is stated as follows.

\begin{theorem}[Hausdorff measure estimate of free boundary]
\label{thm:main}
Let $\Omega\subset\mathbb R^n$ be an open set, $n\ge2$ and let $u$ be a nonnegative viscosity solution of
\begin{equation}
\label{eq:main}
\Delta_\infty u=\lambda u_+^m
\qquad\text{in }\Omega,
\end{equation}
where $\lambda>0$ and $0<m<3$.

For any given $z_0\in\Gamma(u)$, there exists $r_0>0$ such that for any $0<r<r_0$, the set $\Gamma(u)\cap B_r(z_0)$ can be covered by at most $MN$ $(n-1)$-dimensional Lipschitz graphs with Lipschitz constants bounded, hence
\[
\mathcal H^{n-1}\bigl(\Gamma(u)\cap B_r(z_0)\bigr)
\le C r^{n-1},
\]
where the constant $M,N\in\mathbb N$ and $C>0$ depend only on the dimension $n$.

Consequently, $\Gamma(u)$ has locally finite
$(n-1)$-dimensional Hausdorff measure and is
countably $(n-1)$-rectifiable.
\end{theorem}

Furthermore, the interior tangent ball
property yields the sharp lower density
estimate
\[
\liminf_{\rho\to0^+}
\frac{
\mathcal L^n(B_\rho(z)\cap\{u>0\})
}{
\mathcal L^n(B_\rho(z))
}
\ge\frac12
\]
at every $z\in\Gamma(u)$.
The constant $1/2$ is optimal.

\subsection*{Organization of the paper}
Section~2 recalls the comparison principle and the basic notions from geometric measure theory. 
In Section~3, we establish the sharp distance-to-free-boundary regularity (upper growth) estimate and restate the matching optimal-constant non-degeneracy (lower growth) estimate of free boundary points.
Section~4 uncovers the geometric structure underlying the sharp growth estimates. We then yield the sharp lower density estimate for the positive phase and prove a quantitative local bound on the $(n-1)$-dimensional Hausdorff measure of the free boundary, thereby solving the open problem raised in \cite{AraujoLeitaoTeixeira2016}.
Finally, in Section~5, we formulate several questions that arise naturally from our analysis and discuss possible extensions to more general reaction terms.

\section{Preliminary Results}

Let $\Omega\subset\R^n$ be an open set, let $\lambda>0$ be a constant and $0<m<3$. For a function $v$, write $v_+:=\max\{v,0\}$.

\begin{definition}[Viscosity solution]
Let $\Omega\subset\mathbb R^n$ be an open set, $\lambda>0$, and $0<m<3$. Consider the equation
\begin{equation}
\label{eq:viscosity-def}
\Delta_\infty u=\lambda u_+^m
\qquad\text{in }\Omega.
\end{equation}

A function $u\in C(\Omega)$ is called a viscosity
subsolution (respectively, supersolution) of
\eqref{eq:viscosity-def} if, for every
$\varphi\in C^2(\Omega)$ such that $u-\varphi$
attains a local maximum (respectively, minimum)
at some $x_0\in\Omega$, one has
\[
\Delta_\infty\varphi(x_0)
\ge \lambda(u(x_0)_+)^m
\qquad
\bigl(\text{respectively, }\le\bigr).
\]

A function $u\in C(\Omega)$ is called a viscosity
solution of \eqref{eq:viscosity-def} if it is both
a viscosity subsolution and a viscosity supersolution.

Furthermore, given $g\in C(\partial\Omega)$,
a function $u\in C(\overline\Omega)$ is called
a viscosity solution of the corresponding
Dirichlet problem if it solves
\eqref{eq:dirichlet} in the viscosity sense
and satisfies $u=\phi$ on $\partial\Omega$.
\end{definition}

We restate the comparison principle of \cite[Lemma~3.2]{AraujoLeitaoTeixeira2016} here.

\begin{lemma}[Comparison principle]
\label{lem:comparison}
Let $\Omega\subset\R^n$ be a bounded domain, and let
$u,v\in C(\overline{\Omega})$ be nonnegative. Suppose that
\[
\Delta_\infty u\geq \lambda u_+^m,
\qquad
\Delta_\infty v\leq \lambda v_+^m
\qquad\text{in }\Omega
\]
in the viscosity sense. If $u\leq v$ on $\partial\Omega$, then
$u\leq v$ in $\Omega$.
\end{lemma}

Here we reacall some basic notions and concepts from geometric measure theory that would be used later to describe the geometric properties of the free boundary.

We denote by $\mathcal{L}^n$ the $n$-dimensional Lebesgue measure and by $\mathcal{H}^s$ the standard $s$-dimensional Hausdorff measure in $\mathbb{R}^n$.

\begin{definition}[Hausdorff measure]
  Let $E\subset\mathbb R^n$ and $s>0$. For each $\delta>0$, the $s$-dimensional Hausdorff measure of $E$ is defined by
\[
\mathcal H^s(E)
:=
\frac{\omega_s}{2^s}
\lim_{\delta\to0^+}
\inf\left\{
\sum_{i=1}^{\infty}(\operatorname{diam}U_i)^s:
E\subset\displaystyle\bigcup_{i=1}^{\infty}U_i,
\operatorname{diam}U_i\le\delta
\right\},
\]
where $\omega_s:=\dfrac{\pi^{s/2}}{\Gamma(1+s/2)}$ is normalized constant so that it agrees with the usual $s$-dimensional Lebesgue measure when $s$ is an integer.
\end{definition}

\begin{definition}[Lipschitz graph]
Let $1\le k\le n-1$ be an integer, and let $V\subset\mathbb R^n$ be a $k$-dimensional linear subspace. Given a subset $D\subset V$ and a Lipschitz mapping $\mathbf{g}:D\to V^\perp$, the set
\[
G:=\{p+\mathbf{g}(p):p\in D\}
\]
is called a $k$-dimensional Lipschitz graph over $V$.

In particular, when $k=n-1$, we may write $V=e^\perp$ for some
$e\in\mathbb S^{n-1}$. Writing $\mathbf{g}(p)=f(p)e$ for a Lipschitz function
$f:D\to\mathbb R$, we obtain
\[
G=\{p+f(p)e:p\in D\}.
\]
In orthonormal coordinates with $e$ as the last coordinate direction,
the graph can be written as
\[
G
=\bigl\{(x_1,\ldots,x_{n-1},f(x_1,\ldots,x_{n-1}))\in\mathbb R^n
\bigr\}.
\]
\end{definition}

The $k=n-1$ form of Lipschitz graphs will be used in our geometric argument later.

We also recall the following standard property of Hausdorff measures under Lipschitz mappings.

\begin{lemma}[Hausdorff measure properties under Lipschitz maps]
Let $E\subset\mathbb R^n$, $s>0$, and let $F:E\to\mathbb R^n$ be an $L$-Lipschitz mapping. Then
\[
\mathcal H^s(F(E))\le L^s\mathcal H^s(E).
\]

In particular, if $D\subset e^\perp$ and $G=\{p+f(p)e:p\in D\}$ is an
$(n-1)$-dimensional Lipschitz graph, then
\[
\mathcal H^{n-1}(G)
\le \bigl(1+\operatorname{Lip}(f)^2\bigr)^{(n-1)/2}\mathcal H^{n-1}(D),
\]
here $F:D\to\mathbb R^n$ is given by $F(p)=p+f(p)e$, and $\operatorname{Lip}(F) \le\sqrt{1+\operatorname{Lip}(f)^2}$.
\end{lemma}

\section{Matching Sharp Growth Estimates}

We first establish a sharp upper bound for $u$ in terms of the distance to the free boundary.

\begin{theorem}[Optimal-constant sharp upper growth estimate]
\label{thm:le-growth}
Let $u$ be a nonnegative viscosity solution of
\eqref{eq:viscosity-def}. For every $y\in\{u>0\}\cap\Omega$ satisfying
\(
\operatorname{dist}(y,\Gamma(u))<
\operatorname{dist}(y,\partial\Omega),
\)
we have
\[
\begin{aligned}
u(y)
&\le \tau(\lambda,m)\operatorname{dist}(y,\Gamma(u))^{\frac{4}{3-m}}\\
&=
\left(\frac{\lambda(3-m)^4}{64(1+m)}\right)^{\frac{1}{3-m}}
\operatorname{dist}(y,\Gamma(u))^{\frac{4}{3-m}}.
\end{aligned}
\]
In particular, both the exponent and the constant are sharp.
\end{theorem}

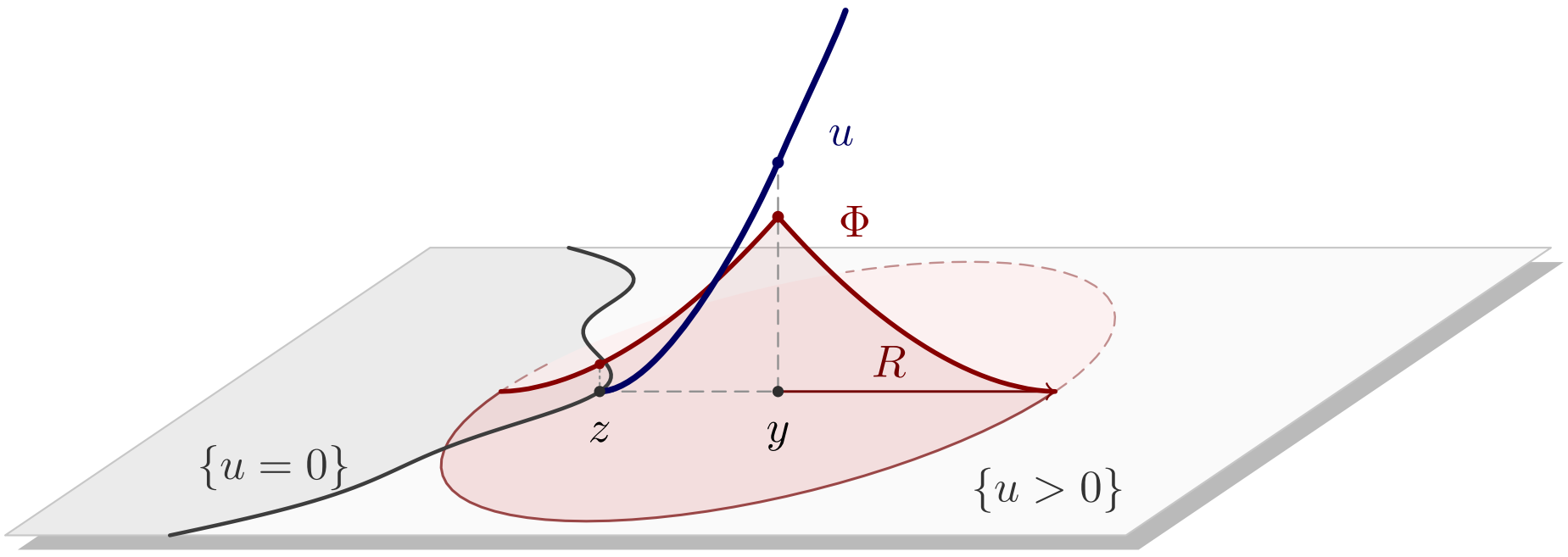
\begin{figure}[htbp]
\centering
\begin{tikzpicture}[
  x={(0.94cm,0cm)}, y={(0.65cm,0.44cm)}, z={(0cm,0.94cm)},
  font=\small, line cap=round, line join=round,
  declare function={
    freeedge(\s)=-1.48-0.27*\s*\s-(0.22+0.06*\s)*(1-cos(170*\s));
  },
  base outline/.style={draw=black!22, line width=.4pt},
  interface/.style={draw=black!76, line width=.85pt},
  solution/.style={draw=blue!39!black, line width=1.35pt},
  barrier/.style={draw=red!53!black, line width=1.05pt},
  guide/.style={draw=black!43, densely dashed, line width=.45pt}
]
\pgfmathsetmacro{\radius}{2.3}
\pgfmathsetmacro{\peak}{1.45}
\pgfmathsetmacro{\profileAtContact}{\peak*(1-abs(freeedge(0))/\radius)^1.8}
\coordinate (center) at (0,0,0);
\coordinate (contact) at ({freeedge(0)},0,0);
\coordinate (barrier at center) at (0,0,\peak);
\coordinate (barrier at contact) at ({freeedge(0)},0,\profileAtContact);
\coordinate (solution at center) at (0,0,1.90);
\coordinate (rim) at (\radius,0,0);

\path[fill=black!27, xshift=2.8pt, yshift=-3.2pt]
  (-4.65,-2.55,0) -- (4.65,-2.55,0) --
  (4.65,2.55,0) -- (-4.65,2.55,0) -- cycle;
\path[fill=black!2] (-4.65,-2.55,0) -- (4.65,-2.55,0) --
  (4.65,2.55,0) -- (-4.65,2.55,0) -- cycle;
\path[fill=black!8] (-4.65,-2.55,0) --
  plot[domain=-2.55:2.55, samples=85, variable=\s]
    ({freeedge(\s)},\s,0)
  -- (-4.65,2.55,0) -- cycle;
\draw[base outline] (-4.65,-2.55,0) -- (4.65,-2.55,0) --
  (4.65,2.55,0) -- (-4.65,2.55,0) -- cycle;

\path[fill=red!8, fill opacity=.58]
  plot[domain=0:360, samples=81, variable=\t]
    ({\radius*cos(\t)},{\radius*sin(\t)},0) -- cycle;
\draw[red!45!black, opacity=.42, densely dashed, line width=.45pt]
  plot[domain=0:113, samples=30, variable=\t]
    ({\radius*cos(\t)},{\radius*sin(\t)},0)
  plot[domain=168:180, samples=8, variable=\t]
    ({\radius*cos(\t)},{\radius*sin(\t)},0);
\path[fill=red!55!black, fill opacity=.08]
  plot[domain=-\radius:\radius, samples=61, variable=\s]
    (\s,0,{\peak*(max(0,1-abs(\s)/\radius))^1.8})
  -- plot[domain=0:-180, samples=41, variable=\t]
    ({\radius*cos(\t)},{\radius*sin(\t)},0) -- cycle;
\draw[red!45!black, opacity=.70, line width=.58pt]
  plot[domain=180:360, samples=41, variable=\t]
    ({\radius*cos(\t)},{\radius*sin(\t)},0);

\draw[interface] plot[domain=-2.55:2.55, samples=85, variable=\s]
  ({freeedge(\s)},\s,0);
\draw[guide] (contact) -- (center);
\draw[guide] (center) -- (solution at center);
\draw[guide, densely dotted] (contact) -- (barrier at contact);
\draw[red!45!black, ->, line width=.55pt]
  (center) -- node[pos=.4,above=2pt,inner sep=1pt] {$R$} (rim);

\draw[barrier] plot[domain=-\radius:\radius, samples=61, variable=\s]
  (\s,0,{\peak*(max(0,1-abs(\s)/\radius))^1.8});
\draw[solution]
  (contact) .. controls (-.94,0,.012) and (-.22,0,1.39) ..
  (solution at center) .. controls (.22,0,2.41) and (.45,0,2.86) ..
  (.56,0,3.16);

\fill[black!82] (contact) circle[radius=1.25pt]
  (center) circle[radius=1.25pt];
\fill[red!53!black] (barrier at contact) circle[radius=1.10pt]
  (barrier at center) circle[radius=1.30pt];
\fill[blue!39!black] (solution at center) circle[radius=1.30pt];

\node[below=3pt] at (contact) {$z$};
\node[below=3pt] at (center) {$y$};
\node[blue!39!black, inner sep=1pt, xshift=14pt, yshift=6pt]
  at (solution at center) {$u$};
\node[red!53!black, inner sep=1pt, xshift=17pt, yshift=-1pt]
  at (barrier at center) {$\Phi$};
\node[black!80] at (-3.25,-1.35,0) {$\{u=0\}$};
\node[black!80] at (3.45,-1.75,0) {$\{u>0\}$};
\end{tikzpicture}
\caption{Radial comparison for the sharp upper growth estimate.}
\label{fig:sharp-upper-comparison}
\end{figure}

\begin{proof}
We argue by contradiction. Figure~\ref{fig:sharp-upper-comparison}
illustrates the radial comparison.

Suppose there exists $y\in\{u>0\}\cap\Omega$ such that
\[
u(y)>\tau\operatorname{dist}(y,\Gamma(u))^{\frac{4}{3-m}}.
\]

By the strict inequality, we can take suitable $R$ satisfying $\operatorname{dist}(y,\Gamma(u))<R<\operatorname{dist}(y,\partial\Omega)$ and $\widetilde{\tau}>\tau$ such that
\[
u(y)>\widetilde{\tau}R^{\frac{4}{3-m}}.
\]

Define inside $\overline{B_R(y)}$ the radical barrier function
\[
\Phi(x):=\widetilde{\tau}(R-|x-y|)^{\frac{4}{3-m}},
\]

Notice that $\Phi$ is not a viscosity subsolution at $y$. Indeed, the constant function $\varphi\equiv\Phi(y)$ touches $\Phi$ from above at $y$, whereas
\[
\Delta_\infty\varphi(y)=0<\lambda\Phi(y)^m.
\]

We therefore apply the comparison principle in a punctured ball.

Since $u(y)>\Phi(y)$ and both $u$ and $\Phi$ are continuous, we may choose $\varepsilon$ sufficiently small such that 
\[
\Phi<u\qquad\text{on }\partial B_\varepsilon(y).
\]
On the other hand,
\[
\Phi=0\le u\qquad\text{on }\partial B_R(y).
\]

Applying the Comparison Principle \ref{lem:comparison}, we can obtain
\[
\Phi\le u
\qquad\text{in }
B_R(y)\setminus\overline{B_\varepsilon(y)}.
\]

Note that there exists $z\in\Gamma(u)$ such that $|z-y|=\operatorname{dist}(y,\Gamma(u))$.

Consequently,
\[
0=u(z)\ge\Phi(z)=
\widetilde{\tau}(R-|z-y|)^{\frac{4}{3-m}}>0,
\]
leading to a contradiction.

Therefore we obtain the upper bound regularity property,
\[
u(y)\le
\tau\operatorname{dist}(y,\Gamma(u))^{\frac{4}{3-m}}.
\]

The sharpness of both the exponent and the constant follows from the
explicit solution $U(x)=\tau(x\cdot e)_+^{\frac{4}{3-m}}$
($e\in\mathbb S^{n-1}$), which attains equality in the estimate.
\end{proof}

\begin{remark}
The improvement-of-flatness and iterative scaling methods or Harnack inequality methods have been widely employed in the studies of regularity estimates to establish optimal growth exponents for this kind of equations; see \cite{AraujoLeitaoTeixeira2016,DiehlTeymurazyan2020,LinLiu2022}, \cite{BezerraDaSilvaNascimentoSa2024,ChenJiang2026}.

For the present model, the direct comparison method gives a simpler proof and a sharper quantitative outcome: the coefficient is explicitly determined and optimal, independently of the $L^\infty$-norm of the solution. Later we will see the great power of the optimal coefficient in the geometric properties reasoning below.

In fact, similar argument can be used to improve the regularity estimate theorem in \cite{ChenJiang2026}, where the constant depends on $\|u\|_{\infty}$. This, however, lead to a natural question that for H\'enon-type equation whether the optimal constant of non-degeneracy estimates
\[
C_{\alpha,m}:=\left(\frac{(3-m)^4}
{(4+\alpha)^3(1+\alpha+m)}
\right)^{\frac1{3-m}}.
\]
can be achieved in the upper bound growth estimate taking account the spatially nonnegative varying weight.

\end{remark}

The non-degeneracy estimate was established by Araújo, Leitão, and Teixeira in \cite[Theorem 6.1]{AraujoLeitaoTeixeira2016} for points in the closure of the positive phase. In particular, the result applies to free boundary points. For our purposes, we restate this estimate with the optimal constant $\tau$, following their barrier argument.

\begin{theorem}[Optimal-constant non-degeneracy at free-boundary points]
\label{thm:ge-growth}
  Let $u\in C(\Omega)$ be a nonnegative viscosity solution of \eqref{eq:main}. Then, for every $z\in\Gamma(u)$ and $0<r<\operatorname{dist}(z,\partial\Omega)$,
\[
\sup_{\partial B_r(z)}u
\ge \tau(\lambda,m)r^{\frac{4}{3-m}}
=
\left(\frac{\lambda(3-m)^4}{64(1+m)}\right)^{\frac{1}{3-m}}
r^{\frac{4}{3-m}}.
\]
The constant $\tau(\lambda,m)$ is optimal.
\end{theorem}

\begin{proof}
We first prove for points within the set $\{u>0\}\cap\Omega$.

Fix $z\in\{u>0\}\cap\Omega$ and
$0<r<\operatorname{dist}(y,\partial\Omega)$, we define
\[
\Psi(x):=\tau|x-z|^{\frac{4}{3-m}}.
\]
Direct computation and by the definition of viscosity solution yield that
$\Delta_\infty\Psi=\lambda\Psi^m$
in the viscosity sense.

We claim that there exist a point $z_r\in\partial B_r(z)$ such that $\Psi \le u$. Otherwise by Comparison Principle \ref{lem:comparison} $\Psi \ge u$ in the entire ball $B_r(z)$, which leads a contradiction
\(
0<u(z)\le\Psi(z)=0.
\)

Therefore we obtain that
\(
\sup_{\partial B_r(z)}u>\tau r^\frac{4}{3-m}.
\)

Then we consider the free boundary points $z\in\Gamma(u)$. Set $\tilde{r}=r-|z_k-z|$, take $z_k\in\{u>0\}$ such that $z_k\to z$ and it holds for every $z_k$ that
\[
\sup_{\partial B_{\tilde{r}}(z_k)}u > \tau \tilde{r}^\frac{4}{3-m} = \tau (r-|z_k-z|)^\frac{4}{3-m}.
\]

Letting $k\to\infty$ and we yield that
\[
\sup_{\partial B_r(z)}u \ge \tau(\lambda,m)r^\frac{4}{3-m}
=
\left(\frac{\lambda(3-m)^4}{64(1+m)}\right)^{\frac{1}{3-m}}
r^{\frac{4}{3-m}}.
\]

The optimality of $\tau$ follows from the same explicit solution $U$, which attains equality in the estimate at $z=0$, showing that $\tau$ is optimal.

\end{proof}

\section{Geometric estimates}
The geometric argument is closely related to that in \cite[Corollary 6.2]{AraujoLeitaoTeixeira2016}, where the regularity property and non-degeneracy estimates are combined to establish the porosity of the free boundary. We have already obtained that upper and lower bound share the same optimal coefficient, which leads to a saturation of the growth inequalities, yielding a uniform interior tangent ball scenery. We state it as a theorem here.

Moreover, a similar geometric idea appears in \cite[Proposition 2.4]{FuscoGelliPisante2012}, where a uniform interior cone condition at each boundary point is used to establish a finite covering of the boundary by Lipschitz graphs. Our argument is quite similar in the sense of structure, but proceeds directly from the uniform interior tangent ball property.

\begin{theorem}[Uniform interior tangent ball]
\label{thm:uniform-tangent}
Let $u$ be a nonnegative viscosity solution of \eqref{eq:main}. For every $z\in\Gamma(u)$ and every $r>0$ such that $\overline{B_{2r}(z)}\subset\Omega$, there exists $y_r\in\partial B_r(z)$ such that
\[
B_r(y_r)\subset\{u>0\},
\qquad
z\in\partial B_r(y_r)\cap\Gamma(u),
\]
and $u(y_r)=\tau(\lambda,m)r^\frac{4}{3-m}$. In fact, $B_r(y_r)$ is an interior ball tangent to the free boundary at $z$.
\end{theorem}

\begin{proof}
By Theorem \ref{thm:ge-growth}, there exists $y_r\in\partial B_r(z)$ and every $r>0$ such that $u(y_r)\ge\tau r^\frac{4}{3-m}$. Applying Theorem \ref{thm:le-growth} at $y_r$, it holds that
\[
\tau r^\frac{4}{3-m}
\le u(y_r)
\le\tau\operatorname{dist}(y_r,\Gamma(u))^\frac{4}{3-m}
\le\tau r^\frac{4}{3-m}.
\]
The inequality chain holds if and only if all the equalities are obtained.

Hence $\operatorname{dist}(y_r,\Gamma(u))=r$, completing the proof.
\end{proof}

The uniform tangent ball property at free boundary points leads to a sharp improvement of the porosity property of the equation.

Araújo, Leitão and Teixeira \cite[Corollary~6.2]{AraujoLeitaoTeixeira2016} established a uniform positive density estimate for the positive phase. More precisely, let $u\in C(\Omega)$ be a viscosity solution of \eqref{eq:main} where $\lambda>0$ and $0\leq m<3$.
Then for every subdomain $\Omega'\Subset\Omega$, there exists a constant $\iota>0$, depending on $\|u\|_{L^\infty(\Omega)}$, $m$ and $\Omega'$, such that
\[
\frac{\mathcal{L}^n(B_r(X_0)\cap\{u>0\})}{r^n}
\geq\iota
\]
for any $X_0\in\partial\{u>0\}\cap\Omega'$
and $0<r\ll1$.

Here we state the improved proposition.

\begin{theorem}[Sharp Lower Density Estimate for the Positive Phase]
Let $u\in C(\Omega)$ be a nonnegative viscosity solution of \eqref{eq:main}. 
Then, for every $z\in\Gamma(u)$,
\[
\liminf_{\rho\to0^+}
\frac{\mathcal L^n(B_\rho(z)\cap\{u>0\})}
     {\mathcal L^n(B_\rho(z))}
\ge\frac12.
\]
Moreover, $\frac12$ here is optimal.
\end{theorem}

\begin{figure}[!htbp]
\centering
\begin{tikzpicture}[
  scale=1.06, font=\small, line cap=round, line join=round,
  tangent ball/.style={draw=red!53!black, line width=.82pt},
  sample ball/.style={draw=blue!39!black, line width=.75pt},
  interface/.style={draw=black!76, line width=.95pt},
  inward/.style={-stealth, draw=blue!39!black, line width=.95pt},
  guide/.style={draw=black!40, dash pattern=on 2.8pt off 2.2pt, line width=.45pt}
]
\def\tangentRadius{1.80}
\def\sampleRadius{.82}
\coordinate (center) at (0,\tangentRadius);

\begin{scope}[rotate around={16:(center)}]
  \coordinate (z) at (0,0);
  \fill[red!7] (center) circle[radius=\tangentRadius];
  \fill[blue!3] (z) circle[radius=\sampleRadius];
  \begin{scope}
    \clip (z) circle[radius=\sampleRadius];
    \fill[red!23] (center) circle[radius=\tangentRadius];
  \end{scope}

  \draw[guide] (-\sampleRadius,0) -- (\sampleRadius,0);
  \draw[tangent ball] (center) circle[radius=\tangentRadius];
  \draw[sample ball] (z) circle[radius=\sampleRadius];
  \draw[interface]
    (-3.30,.06) .. controls (-2.75,.10) and (-2.17,-.53) .. (-1.56,-.50)
    .. controls (-.94,-.46) and (-.54,0) .. (z)
    .. controls (.51,0) and (.99,-.49) .. (1.44,-.46)
    .. controls (1.80,-.44) and (2.01,-.28) .. (2.31,-.33);

  \draw[inward] (z) -- node[pos=.70, right=4.5pt, text=blue!39!black]
    {$r\nu_z$} (center);
  \fill[black!82] (z) circle[radius=1.25pt] node[below right=2pt] {$z$};
\end{scope}
\fill[black!82] (center) circle[radius=1.25pt];
\node[red!53!black] at (0,2.65) {$B_r(z+r\nu_z)$};
\node[blue!39!black, anchor=west] at (1.24,-.56) {$B_\rho(z)$};
\node[black!76, anchor=north] at (-1.79,-1.03) {$\{u=0\}$};
\node[black!76] at (2.37,.84) {$\{u>0\}$};
\end{tikzpicture}
\caption{The overlap $B_\rho(z)\cap B_r(z+r\nu_z)$ used in the lower density estimate.}
\label{fig:tangent-ball-density}
\end{figure}
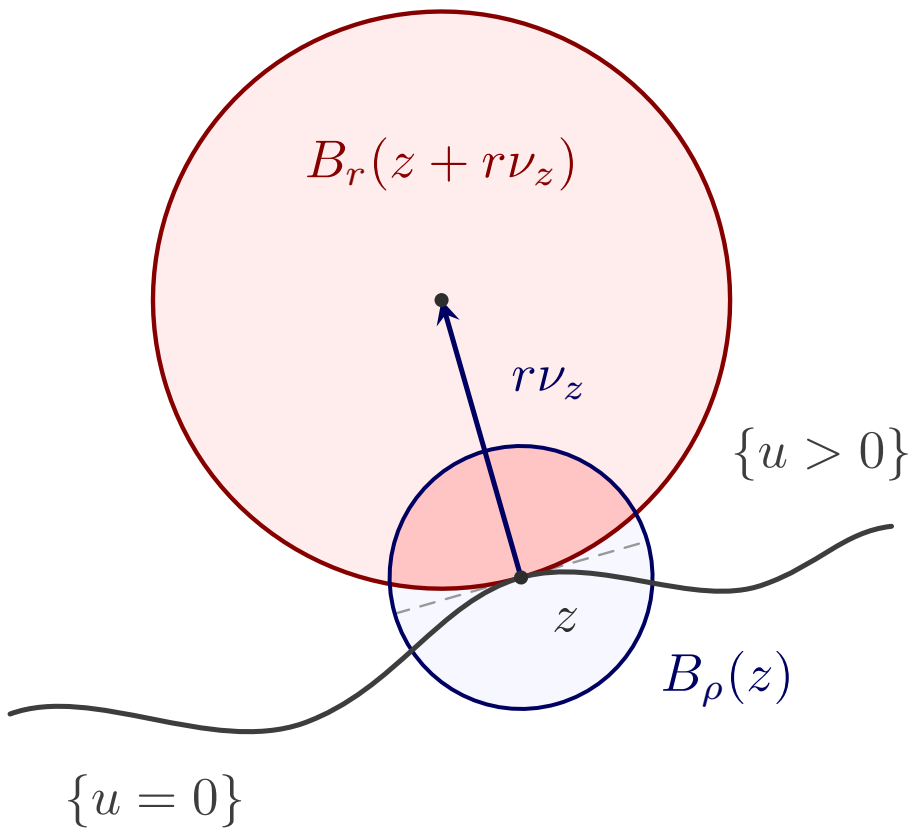

\begin{proof}
For any fixed free boundary point $z\in\Gamma(u)$ and $0<\rho<r$, consider $x=z+\rho y$.

We can describe the geometric relation using the inequalities below,
\[
x\in B_\rho(z)\iff |y|<1,
\]
\[
\begin{aligned}
x\in B_r(z+r\nu_z)
&\iff |\rho y-r\nu_z|^2<r^2\\
&\iff y\cdot\nu_z>\frac{\rho}{2r}|y|^2.
\end{aligned}
\]

Define
\(
D_\rho=
\left\{
y\in B_1(0):
y\cdot\nu_z>\frac{\rho}{2r}|y|^2
\right\}
\)
and
\(
D_0=B_1(0)\cap\{y:y\cdot\nu_z>0\}.
\)

Consider
\[
\begin{aligned}
\frac{\mathcal L^n(B_\rho(z)\cap\{u>0\})}
{\mathcal L^n(B_\rho(z))}
&\ge
\frac{\mathcal L^n\bigl(B_\rho(z)\cap B_r(z+r\nu_z)\bigr)}
{\omega_n\rho^n}\\
&=\frac1{\omega_n}\int_{D_\rho}dy,
\end{aligned}
\]

Note that as $\rho\to 0^+$,
\[
\lim_{\rho\to0^+}\frac1{\omega_n}\int_{D_\rho}dy
=
\frac{\mathcal L^n(D_0)}{\omega_n}
=\frac12.
\]

Therefore we obtain the density estimate,
\[
\liminf_{\rho\to0^+}
\frac{\mathcal L^n(B_\rho(z)\cap\{u>0\})}
{\mathcal L^n(B_\rho(z))}
\ge\frac12.
\]
\end{proof}

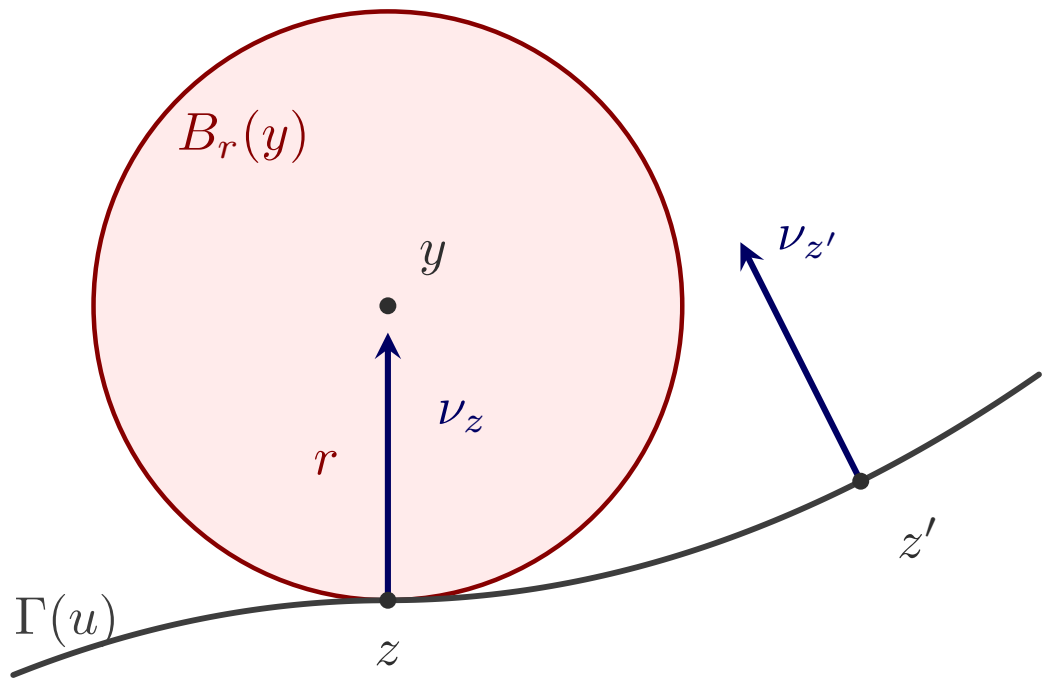
\begin{figure}[htbp]
\centering
\begin{tikzpicture}[
  scale=1.18, font=\small, line cap=round, line join=round,
  ball/.style={draw=red!53!black, fill=red!8, line width=.85pt},
  interface/.style={draw=black!76, line width=1.10pt},
  inward/.style={-stealth, draw=blue!39!black, line width=1.15pt}
]
\coordinate (z) at (0,0);
\coordinate (y) at (0,1.65);
\coordinate (zp) at (2.65,{0.095*2.65^2});

\path[ball] (y) circle[radius=1.65];
\draw[interface, smooth, domain=-2.1:3.65, samples=61]
  plot (\x,{0.095*\x^2});

\def\normalLength{1.50}
\draw[inward] (z) -- ++(90:\normalLength);
\draw[inward] (zp) -- ++({90+atan(2*0.095*2.65)}:\normalLength)
  coordinate (zp-tip);

\fill[black!82] (z) circle[radius=1.35pt] node[below=4pt] {$z$};
\fill[black!82] (y) circle[radius=1.35pt] node[above right=2pt] {$y$};
\fill[black!82] (zp) circle[radius=1.35pt] node[below right=3pt] {$z'$};
\node[red!53!black, anchor=east] at (-0.16,0.77) {$r$};
\node[blue!39!black, anchor=west] at (0.16,1.03) {$\nu_z$};
\node[blue!39!black, right=3pt] at (zp-tip) {$\nu_{z'}$};
\node[red!53!black] at (-0.81,2.60) {$B_r(y)$};
\node[black!76, anchor=north] at (-1.80,0.15) {$\Gamma(u)$};
\end{tikzpicture}
\caption{The interior tangent ball $B_r(y)$ at $z$ and the inward
directions $\nu_z$ and $\nu_{z'}$ at two free boundary points.}
\label{fig:interior-tangent-ball}
\end{figure}

The uniform interior tangent ball property established above also provides quantitative control of the geometry of the free boundary itself. Indeed, we can show that the free boundary can be covered by finitely many Lipschitz graphs. This leads to the desired quantitative $(n-1)$-dimensional Hausdorff measure estimate.

We now prove our main theorem.

\begin{proof}[Proof of Theorem \ref{thm:main}]
Fix $z_0\in\Gamma(u)$ and take suitable $r_0>0$.

For $0<r<r_0$, Theorem \ref{thm:uniform-tangent} gives the interior tangent ball relation. 

For any $z,z'\in\Gamma(u)\cap B_r(z_0)$, we obtain
\[
\nu_z\cdot(z'-z)\le\frac{|z'-z|^2}{2r},
\]
and
\[
\nu_{z'}\cdot(z-z') \le\frac{|z'-z|^2}{2r}.
\]

We classify the free boundary points into finitely many classes.

First, consider finitely many directions
\(
e_1,\ldots,e_M\in\mathbb S^{n-1},
\)
where $M<\infty$ depends only on $n$, such
that for any $\nu_z$, there exists an $e_i$
satisfying
\[
|\nu_z-e_i|<\frac14.
\]

Meanwhile, we can choose $N$ points $x_j$ such that
\(
B_r(z_0)\subset\bigcup_{j=1}^N B_{r/4}(x_j),
\)
where $N$ depends only on $n$.
For $1\le i\le M$ and $1\le j\le N$, set
\[
\Gamma_{ij}:=\bigl\{z\in\Gamma(u)\cap B_r(z_0)\cap B_{r/4}(x_j):
|\nu_z-e_i|<\tfrac14\bigr\}.
\]
These classes cover $\Gamma(u)\cap B_r(z_0)$.

Consider a nonempty class $\Gamma_{ij}$.
For $z,z'\in\Gamma_{ij}$, with $e_i$, we have
\[
|\nu_z-e_i|<\frac14,\qquad
|\nu_{z'}-e_i|<\frac14,
\qquad |z-z'|\le\frac r2.
\]

Using the two tangency inequalities above,
we have
\[
\begin{aligned}
|e_i\cdot(z-z')|
&\le \frac14|z-z'|
+\frac{|z-z'|^2}{2r}\\
&\le \frac12|z-z'|.
\end{aligned}
\]

Let $P_i$ be the orthogonal projection
onto $e_i^\perp$, which is
\begin{align*}
P_i:\mathbb R^n&\longrightarrow e_i^\perp,\\
z&\longmapsto z-(z\cdot e_i)e_i.
\end{align*}
In particular, $P_i$ is injective defined on $\Gamma_{ij}$ and 
\begin{align*}
|z-z'|^2
&=|P_i(z-z')|^2
  +|e_i\cdot(z-z')|^2\\
&\le |P_i z-P_i z'|^2
  +\frac14|z-z'|^2.
\end{align*}

Set $D_{ij}:=P_i(\Gamma_{ij})\subset e_i^\perp.$ and define ``height function''
\[
f_{ij}:D_{ij}\to\mathbb R,
\qquad
f_{ij}(P_i z):=z\cdot e_i,
\quad z\in\Gamma_{ij}.
\]

For any given $z,z'\in\Gamma_{ij}$, we have
\begin{align*}
|f_{ij}(P_i z)-f_{ij}(P_i z')|
&=|e_i\cdot(z-z')|\\
&\le\frac12|z-z'|\\
&\le\frac1{\sqrt3}|P_i z-P_i z'|.
\end{align*}

Then we can define the graph map
\begin{align*}
F_{ij}:D_{ij}&\longrightarrow\mathbb R^n,\\
p&\longmapsto p+f_{ij}(p)e_i.
\end{align*}

Note that $\Gamma_{ij}=F_{ij}(D_{ij})$ and 
\[
\operatorname{Lip}(F_{ij})
\le\sqrt{1+\operatorname{Lip}(f_{ij})^2}
\le\frac2{\sqrt3},
\]
yielding that $\Gamma_{ij}$ is an $(n-1)$-dimensional Lipschitz graph with the graph map $F_{ij}$.

Then, applying the Hausdorff measure estimate for a Lipschitz graph,
\begin{align*}
\mathcal H^{n-1}(\Gamma_{ij})
&=\mathcal H^{n-1}(F_{ij}(D_{ij}))\\
&\le\left(\frac2{\sqrt3}\right)^{n-1}
\mathcal H^{n-1}(D_{ij})\\
&\le\left(\frac2{\sqrt3}\right)^{n-1}
\omega_{n-1}\left(\frac r4\right)^{n-1}=\omega_{n-1}
\left(\frac r{2\sqrt3}\right)^{n-1},
\end{align*}
where $\omega_{n-1}$ is the volume of the
unit ball in $\mathbb R^{n-1}$.

Summing over the finitely many classes,
we obtain
\[
\begin{aligned}
\mathcal H^{n-1}
  \bigl(\Gamma(u)\cap B_r(z_0)\bigr)
  &\le\sum_{i=1}^M\sum_{j=1}^N
  \mathcal H^{n-1}(\Gamma_{ij})\\
  &\le MN\,\omega_{n-1}
  \left(\frac r{2\sqrt3}\right)^{n-1}=:Cr^{n-1},
\end{aligned}
\]
where $C>0$ depends only on $n$.

This gives the finite Lipschitz-graph covering by at most $MN$ graphs
with Lipschitz constants bounded by $2/\sqrt{3}$.
The stated local finiteness and countable
rectifiability follow.
\end{proof}

\section{Final comments and possible generalizations}

This work is devoted to the nonnegative viscosity solutions and the corresponding geometric properties of the free boundary.

Ara\'ujo, Leit\~ao and Teixeira~ established the local finiteness of the $(n-\sigma_0)$-dimensional Hausdorff measure of the free boundary for some $\sigma_0>0$; see \cite[Remark~6.3]{AraujoLeitaoTeixeira2016}.
They conjectured there that one could take $\sigma_0=1$ and left this as an open problem.

We give an affirmative answer to this question. The key observation is that the sharp upper growth estimate and the non-degeneracy estimate share the same optimal coefficient \(\tau\), which plays a key role in the geometric argument.

Hence, a natural question is the optimal Hausdorff measure constant. We already obtain the quantitative estimate of the Hausdorff measure
\[
\mathcal H^{n-1}\bigl(\Gamma(u)\cap B_r(z_0)\bigr)
\le C r^{n-1},
\]
while the constant $C$ is quite rough and we want to know the optimal constant, defined as
\[
C_n^{\mathrm{opt}}
:=
\sup
\frac{
\mathcal H^{n-1}
\bigl(\Gamma(u)\cap B_r(z)\bigr)
}{r^{n-1}},
\]
and by Theorem \ref{thm:main}, $C_n^{\mathrm{opt}}<\infty$. Also, the constant may have something to do with the shape of the free boundary.

Another natural extension is consider weight functions spatially varying, which is
\[
\Delta_\infty u=\lambda(x)u_+^m,
\qquad 0<m<3,
\]
provided $\lambda(x)$ with suitable conditions.

For exmaple, like existing literature considering $\lambda\in C(\Omega)$ and $\inf_\Omega \lambda>0$, we can ask which properties still hold under these assumptions.

To be more precise, for fixed $z\in\Gamma(u)$ and sufficiently small $r>0$, define
\[
\lambda_0:=
\min_{\overline{B_{3r}(z)}}\lambda,
\qquad
\Lambda:=
\max_{\overline{B_{3r}(z)}}\lambda,
\]
thus $0<\lambda_0\le\lambda(x)\le\Lambda$ inside $\overline{B_{3r}(z)}$ and
\[
\lambda_0u_+^m\le\Delta_\infty u\le\Lambda u_+^m.
\]

Furthermore, applying methods used in establishing the estimates for the solutions and taking $u(y_r)=\sup_{B_r(z)}u$, it holds that
\[
\left(\frac{\lambda_0(3-m)^4}{64(1+m)}\right)^{\frac1{3-m}}r^{\frac4{3-m}}
\le u(y_r)
\le
\left(\frac{\Lambda(3-m)^4}{64(1+m)}\right)^{\frac1{3-m}}\operatorname{dist}(y_r,\Gamma(u))^{\frac4{3-m}}.
\]

Thus,
\[
\operatorname{dist}(y_r,\Gamma(u))
\ge
\left(\frac{\lambda_0}{\Lambda}\right)^{1/4}r:=\theta_z(r)\cdot r.
\]

In other words, this indicates that
\[
B_{\theta_z r}(y_r)
\subset\{u>0\},
\qquad y_r\in\partial B_r(z).
\]

As we suppose $\lambda(x)\in C(\Omega)$ and $\inf_\Omega \lambda>0$ here, the fact that $\lim_{r\to0^+}\theta_z(r)=1$ leads to a ``almost tangent ball'' property, indicating similar structure for strong absorption equations with varying weight and it worth further research.

Finally, another natural direction is to consider equations with more general reaction terms, beginning with autonomous equations and extending to spatially and gradient-dependent nonlinearities,
\[
\Delta_\infty u=f(u),\qquad
\Delta_\infty u=f(x,u),\qquad
\Delta_\infty u=f(x,u,Du).
\]
where $f$ can have some good properties. 
For example, related regularity, non-degeneracy, and porosity estimates have been investigated by Diehl and Teymurazyan, see \cite{DiehlTeymurazyan2020}. These questions merit further investigation.

\section*{Statement on the use of AI}

The authors identified the research direction and provided the relevant mathematical background. OpenAI's ChatGPT and Codex assisted in exploring proof strategies and generated an initial proof. ChatGPT was also used for language editing. The authors checked the references, reviewed and finalized the manuscript, and take full responsibility for its mathematical correctness and content.

\ifreferencesused
\bibliographystyle{amsplain-doi}
\bibliography{refs}
\else
\section*{References}
\fi

\end{document}